\documentclass[a4paper,11pt]{article}
\usepackage{amsmath,amscd,amsfonts,amssymb,epsf,latexsym}
\usepackage[all]{xy}
\usepackage{graphicx}

\makeatletter

\long\def\@makefnt#1{\parindent 1em\noindent
            \hb@xt@1.8em{\hss\@textsuperscript{}}#1}
\long\def\@ftntext#1{\insert\footins{%
    \reset@font\footnotesize
    \interlinepenalty\interfootnotelinepenalty
    \splittopskip\footnotesep
    \splitmaxdepth \dp\strutbox \floatingpenalty \@MM
    \hsize\columnwidth \@parboxrestore
    \color@begingroup
      \@makefnt{%
        \rule\z@\footnotesep\ignorespaces#1\@finalstrut\strutbox}%
    \color@endgroup}}%
\def\subjclass#1{%
  \@ftntext{2010 {\itshape Mathematics Subject Classification.}\enspace #1.}}
\def\keywords#1{%
  \@ftntext{{\itshape Key words and phrases.}\enspace #1.}}
\makeatother

\def\moins{\raise 1pt\hbox{{$\scriptstyle -$}}}
\def\plus{\raise 1pt\hbox{{$\scriptstyle +$}} }

\newtheorem{theorem}{Theorem}
\newtheorem{proposition}[theorem]{Proposition}
\newtheorem{lemma}[theorem]{Lemma}
\newtheorem{corollary}[theorem]{Corollary}
\newtheorem{remark}[theorem]{Remark}

\newtheorem{example}[theorem]{Example}
\newtheorem{note}[theorem]{Note}

\def\proof{\noindent{\bf Proof.\ }}

\def\qed{~\hbox{$\Box$}}

\def\rank{\mathop{\rm rank}}
\def\dim{\mathop{\rm dim}}

\def\Hom{\mathop{\rm Hom}}
\def\Ext{\mathop{\rm Ext}}

\def\Ker{\mathop{\rm Ker}}

\def\Z{{\mathbb Z}}

\def\rank{\mathop{\mathrm{rank}}}
\def\Ker{\mathop{\mathrm{Ker}}}

\let\la\lambda

\def\N{{\mathbb N}}

\def\Z{{\mathbb Z}}

\def\cC{{\mathcal C}}

\def\Hom{\mathop{\rm Hom}}

\def\moins{\raise 1pt\hbox{{$\scriptstyle -$}}}
\def\plus{\raise 1pt\hbox{{$\scriptstyle +$}} }

\def\phi{\varphi}

\catcode`\@=11
\def\binrel@#1{\setbox\z@\hbox{\thinmuskip0mu
 \medmuskip\m@ne mu\thickmuskip\@ne mu$#1\m@th$}%
 \setbox\@ne\hbox{\thinmuskip0mu\medmuskip\m@ne mu\thickmuskip
 \@ne mu${}#1{}\m@th$}%
 \setbox\tw@\hbox{\hskip\wd\@ne\hskip-\wd\z@}}

\def\overset#1\to#2{\binrel@{#2}\ifdim\wd\tw@<\z@
 \mathbin{\mathop{\kern\z@#2}\limits^{#1}}\else\ifdim\wd\tw@>\z@
 \mathrel{\mathop{\kern\z@#2}\limits^{#1}}\else
 {\mathop{\kern\z@#2}\limits^{#1}}{}\fi\fi}

\def\underset#1\to#2{\binrel@{#2}\ifdim\wd\tw@<\z@
 \mathbin{\mathop{\kern\z@#2}\limits_{#1}}\else\ifdim\wd\tw@>\z@
 \mathrel{\mathop{\kern\z@#2}\limits_{#1}}\else
 {\mathop{\kern\z@#2}\limits_{#1}}{}\fi\fi}
 \catcode`\@=12

\def\POZ#1,#2)#3.{\put(#1,#2){\line(1,0){#3}}}
\def\PION#1,#2)#3.{\put(#1,#2){\line(0,1){#3}}}

\def\proof{\noindent{\bf Proof.\ }}
\def\qed{~\hbox{$\Box$}}

\def\binrel@#1{\setbox\z@\hbox{\thinmuskip0mu
 \medmuskip\m@ne mu\thickmuskip\@ne mu$#1\m@th$}%
 \setbox\@ne\hbox{\thinmuskip0mu\medmuskip\m@ne mu\thickmuskip
 \@ne mu${}#1{}\m@th$}%
 \setbox\tw@\hbox{\hskip\wd\@ne\hskip-\wd\z@}}

\def\overset#1\to#2{\binrel@{#2}\ifdim\wd\tw@<\z@
 \mathbin{\mathop{\kern\z@#2}\limits^{#1}}\else\ifdim\wd\tw@>\z@
 \mathrel{\mathop{\kern\z@#2}\limits^{#1}}\else
 {\mathop{\kern\z@#2}\limits^{#1}}{}\fi\fi}

\def\underset#1\to#2{\binrel@{#2}\ifdim\wd\tw@<\z@
 \mathbin{\mathop{\kern\z@#2}\limits_{#1}}\else\ifdim\wd\tw@>\z@
 \mathrel{\mathop{\kern\z@#2}\limits_{#1}}\else
 {\mathop{\kern\z@#2}\limits_{#1}}{}\fi\fi}
 \catcode`\@=12

\begin{document}

\title{\bf On a certain family of $U(\frak b)$-modules}

\author{Piotr Pragacz\thanks{Supported by National Science Center (NCN)
grant no. 2014/13/B/ST1/00133}\\
\small Institute of Mathematics, Polish Academy of Sciences\\
\small \'Sniadeckich 8, 00-656 Warszawa, Poland\\
\small P.Pragacz@impan.pl}



\date{}

\maketitle



\rightline{\it \small To IMPANGA on the occasion of her 15th birthday}

\begin{abstract}
We report on results of Kra\'skiewicz and the author, and Watanabe on
KP modules materializing Schubert polynomials, and filtrations having KP modules as their subquotients.
We discuss applications of KP filtrations and ample KP bundles to positivity due to Fulton and Watanabe
respectively.

\end{abstract}
\section{Introduction}
In the present article, we survey a recent chapter of representation theory of $U(\frak b)$-modules: the theory of KP modules.

There is a really famous family of $U(\frak b)$-modules: these are Demazure modules (\cite{D1}, \cite{D2},
\cite{J}, \cite{A}, \cite{vK}). They are given by the spaces of sections of line bundles on the Schubert varieties in flag manifolds. Invented in the 1970s by Demazure, the theory was developed in the 1980s by Joseph, Andersen, Polo, Mathieu, van der Kallen et al.

The origins of KP modules, though also related to Schubert varieties, are different. In the beginning of 1980s, Lascoux and Sch\"utzenberger \cite{LS1}
discovered Schubert polynomials -- certain polynomial lifts of cohomology classes of Schubert varieties in flag
manifolds (see \cite{BGG}). Schubert polynomials are described in Sect.~\ref{schub}.
It was a conjecture of Lascoux (Oberwolfach, June 1983) that there should exist a functorial
version of this construction (similarly as to Schur functions there correspond Schur functors,
cf. Sect.~\ref{schur}). The Lascoux conjecture was solved affirmatively by Kra\'skiewicz and the author \cite{KP0, KP}.
The so-obtained modules were called {\it Kra\'skiewicz-Pragacz modules}, in short {\it KP modules} by Watanabe, who is the author
of further developments of the theory of KP modules and KP filtrations in the spirit of the highest weight categories \cite{CPS}.
His work \cite{W1, W2, W3} is surveyed in Sect.~\ref{wat}. The KP modules and related modules, e.g. Schur flagged modules, are discussed in Sect. \ref{las}.

The last two sections are devoted to study positivity.

In Section \ref{plet}, we discuss a recent result of Watanabe, showing that a Schur function specialized with the monomials $x_1^{\alpha_1}x_2^{\alpha_2}\cdots$
of a Schubert polynomial is a nonnegative combination of Schubert polynomials in $x_1, x_2,\ldots$. The method relies on KP filtrations.

For a monomial $x^\alpha=x_1^{\alpha_1}x_2^{\alpha_2}\cdots$, set $l(x^{\alpha}):=\alpha_1 x_1 +\alpha_2 x_2 + \cdots$.
In Section \ref{ample}, we discuss a result of Fulton from the 1990s, asserting that a Schur function specialized with the expressions $l(x^{\alpha})$
associated to the monomials $x^{\alpha}$ of a Schubert polynomial is a nonnegative combination of Schubert polynomials 
in $x_1, x_2,\ldots$. The method relies on ample vector bundles.

\smallskip

Perhaps a couple of words about comparison of Demazure and KP modules is in order.
KP modules are in some way similar to Demazure modules (of type (A)), the modules generated
by an extremal vector in an irreducible representation of $\frak gl_n$: they are both cyclic
$U(\frak b)$-modules parametrized by the weights of the generators, and if the permutation is vexillary,
then the KP module coincides with Demazure module with the same weight of the generator. If a permutation
is not vexillary, then there exists a strict surjection from the KP module to the Demazure module
of the corresponding lowest weight.

\smallskip

This article is a written account of the talk given by the author during the conference
``IMPANGA 15'' at the Banach Center in B\c edlewo in April 2015.

\section{Schur functors}\label{schur}

Throughout this paper, let $K$ be a field of characteristic zero.

In his dissertation \cite{Sch} (Berlin, 1901), Schur gave a classification
of irreducible polynomial representations of the full linear group $GL_n(K)$, i.e.
homomorphisms
$$
GL_n(K)\to GL_N(K)
$$
sending an $n\times n$-matrix $X$ to an $N\times N$-matrix $[P_{ij}(X)]$, where $P_{ij}$ is a polynomial
in the entries of $X$. Let $\Sigma_n$ denotes the symmetric group of all bijections of $\{1,\ldots,n\}$. Consider the following two actions on
$E^{\otimes n}$, where $E$ is a finite dimensional vector space over $K$:
\begin{itemize}
\item of the symmetric group $\Sigma_n$ via permutations of the factors;
\item the diagonal action of $GL(E)$.
\end{itemize}

Irreductible representations $S^\la$ of the symmetric group $\Sigma_n$
are labeled by partitions of $n$, see, e.g. \cite{Bo}. By a {\it partition} of $n$, we
mean a sequence
$$
\la=(\la_1\ge \cdots \ge \la_k \ge 0)
$$
of integers such that $|\la|=\la_1+ \cdots+\la_k=n$. A partition is often represented
graphically by its diagram with $\la_i$ boxes in the $i$th row (cf. \cite{Mcd2}).

\begin{example} \rm The diagram of the partition $(8,7,4,2)$ is
$$\unitlength=10pt
\begin{picture}(7,5)(0,0)
\POZ0,0)2.
\POZ0,1)4.
\POZ0,2)7.
\POZ0,3)8.
\POZ0,4)8.
\PION0,0)4.
\PION1,0)4.
\PION2,0)4.
\PION3,1)3.
\PION4,1)3.
\PION5,2)2.
\PION6,2)2.
\PION7,2)2.
\PION8,3)1.
\end{picture}
$$
\end{example}
Let $R$ be a commutative ring at let $E$ be an $R$-module.
We define the {\it Schur module} associated with a partition $\la$
as follows:
$$
V_\la(E):={\Hom}_{\Z[\Sigma_n]}(S^\la, E^{\otimes n})\,,
$$
where $\Z[\Sigma_n]$ denotes the group ring of $\Sigma_n$ with integer coefficients (cf. \cite{Niel}).

In fact, $V_\la(-)$ is a {\it functor}: if $E, F$ are $R$-modules
over a commutative ring $R$, and
$f: E \to F$ is an $R$-homomorphism, then $f$ induces an $R$-homomorphism
$V_\la(E) \to  V_\la(F)$. In this way, we get all irreducible polynomial representations
of $GL_n(K)$.

Let us label the boxes of the diagram of a partition of $n$ with numbers $1,\ldots,n$,
and call such an object a {\it tableau}. For example, a tableau for the diagram of the partition
 $(8,7,4,2)$ is
$$\unitlength=10pt
\begin{picture}(7,5)(0,0)
\POZ0,0)2.
\POZ0,1)4.
\POZ0,2)7.
\POZ0,3)8.
\POZ0,4)8.
\PION0,0)4.
\PION1,0)4.
\PION2,0)4.
\PION3,1)3.
\PION4,1)3.
\PION5,2)2.
\PION6,2)2.
\PION7,2)2.
\PION8,3)1.
\put(0.5,0.5){\makebox(0,0){$\scriptstyle 16$}}
\put(0.5,1.5){\makebox(0,0){$\scriptstyle 7$}}
\put(0.5,2.5){\makebox(0,0){$\scriptstyle 11$}}
\put(0.5,3.5){\makebox(0,0){$\scriptstyle 1$}}
\put(1.5,0.5){\makebox(0,0){$\scriptstyle 2$}}
\put(1.5,1.5){\makebox(0,0){$\scriptstyle 20$}}
\put(1.5,2.5){\makebox(0,0){$\scriptstyle 8$}}
\put(1.5,3.5){\makebox(0,0){$\scriptstyle 15$}}
\put(2.5,1.5){\makebox(0,0){$\scriptstyle 12$}}
\put(2.5,2.5){\makebox(0,0){$\scriptstyle 18$}}
\put(2.5,3.5){\makebox(0,0){$\scriptstyle 19$}}
\put(3.5,1.5){\makebox(0,0){$\scriptstyle 16$}}
\put(3.5,2.5){\makebox(0,0){$\scriptstyle 9$}}
\put(3.5,3.5){\makebox(0,0){$\scriptstyle 3$}}
\put(4.5,2.5){\makebox(0,0){$\scriptstyle 6$}}
\put(4.5,3.5){\makebox(0,0){$\scriptstyle 10$}}
\put(5.5,2.5){\makebox(0,0){$\scriptstyle 17$}}
\put(5.5,3.5){\makebox(0,0){$\scriptstyle 5$}}
\put(6.5,2.5){\makebox(0,0){$\scriptstyle 4$}}
\put(6.5,3.5){\makebox(0,0){$\scriptstyle 21$}}
\put(7.5,3.5){\makebox(0,0){$\scriptstyle 13$}}
\end{picture}$$

Consider the following two elements in $\Z[\Sigma_n]$ associated with a tableau:
\begin{itemize}
\item  $P$:= sum of $w\in \Sigma_n$ preserving the rows of the tableau\,;
\item  $N$:= sum of $w\in \Sigma_n$ with their signs, preserving the columns of the tableau\,.
\end{itemize}

The following element
$$
e(\la):=N \circ P
$$
is called a {\it Young symmetrizer}.
For more on Young symmetrizers, see \cite {Bo}.
We have an alternative presentation of Schur module:
$$
V_\la(E)=e(\la) E^{\otimes n}\,.
$$
There is still another way of getting Schur modules $V_\lambda(E)$ (see \cite{L}, \cite{ABW}) as the images of natural homomorphisms 
between the tensor products of symmetric and exterior powers of a module:
\begin{equation}\label{sw}
S_{\lambda_1}(E) \otimes \cdots \otimes S_{\lambda_k}(E) \to E^{\otimes |\lambda|} \to \wedge^{\mu_1}(E)\otimes \cdots \otimes \wedge^{\mu_l}(E)\,.
\end{equation}
(here $\mu=\lambda^\sim$ is the conjugate partition of $\lambda$, see \cite{Mcd2}, where a different notation is used; the first map is the diagonalization in the symmetric algebra, and the second map is the multiplication in the exterior algebra), and
$$
\wedge^{\mu_1}(E)\otimes \cdots \otimes \wedge^{\mu_l}(E) \to E^{\otimes |\lambda|} \to S_{\lambda_1}(E) \otimes \cdots \otimes S_{\lambda_k}(E) 
$$
(here the first map is the diagonalization in the exterior algebra, and the second map is multiplication in the symmetric algebra).

\begin{example} \rm We have $V_{(n)}(E)=S_n(E)$
and $V_{(1^n)}(E)=\wedge^n(E)$, the $n$th symmetric and exterior power of $E$.
\end{example}

Let $T$ be the subgroup of diagonal matrices in $GL_n(K)$:
\begin{equation}\label{T}
\left(\begin{matrix}
\scriptstyle x_1&&&&\\
&\scriptstyle x_2&&\scriptstyle 0&\\
&&\scriptstyle x_3&&\\
&\scriptstyle0&&\ddots&
\end{matrix}
\right)\,.
\end{equation}
Let $E$ be a finite dimensional vector space over $K$.
Consider the action of $T$ on the Schur module $V_\la(E)$ associated with a partition $\lambda=(\lambda_1,\ldots,\lambda_k)$, induced from the action of $GL_n(K)$
via restriction.

\begin{theorem}[\rm Main result of Schur's Thesis]
The trace of the action of $T$ on $V_\la(E)$ is equal to the Schur function:
$$
s_\la(x)=
\det \left( s_{\la_p-p+q}(x_1,\ldots ,x_n)
\right)_{1\le p,q\le k}\,,
$$
where $s_i(x_1,\ldots,x_n)$, $i\in \Z$, is the complete symmetric function of degree $i$ if $i\ge 0$,
and zero otherwise.
\end{theorem}
For more on Schur functions, see \cite{Mcd2}, \cite{Lbook}.

Among the most important formulas in the theory of symmetric functions, they are
the following Cauchy formulas:
$$
\prod_{i,j}(1-x_i y_j)^{-1}=\sum_\lambda s_\lambda(x) s_\lambda(y) \ \ \ \ \ \hbox{and} \ \ \ \ \
\prod_{i,j}(1+x_i y_j)=\sum_\lambda s_\lambda(x) s_{\lambda^{\sim}}(y)\,.
$$
Schur functors were brought to the attention of algebraists in \cite{L} together with
the materializations\footnote{We say that a formula on the level of representations ``materializes'' a polynomial formula if the latter
is the character of the former.} of the Cauchy formulas: for free $R$-modules $E,F$ one has:
$$
S_n(E\otimes F)=\oplus_{|\lambda|=n} V_\lambda(E)\otimes V_{\lambda}(F)\,,
$$
and
$$
\wedge^n(E\otimes F)=\oplus_{|\lambda|=n} V_\lambda(E)\otimes V_{\lambda^{\sim}}(F)\,.
$$
Whereas the former formula gives the irreducible $GL(E)\times GL(F)$-representa\-tions of the
decomposition of the space of functions on the space of $\dim(E)\times \dim(F)$-matrices,
the importance of the latter comes from the fact that it describes
the {\it Koszul syzygies} of the ideal generated by the entries
of a generic $\dim(E)\times \dim(F)$-matrix. Then, using suitable derived functors
(following a method introduced in Kempf's 1971 Thesis), this 
allows one to describe syzygies of determinantal ideals \cite{L}.
In fact, there is a natural extension of Schur functors to Schur complexes  
with many applications (see \cite{Niel}, \cite{ABW}).

\section{Schubert polynomials}\label{schub}

In the present article, by a {\it permutation} $w=w(1),w(2),\ldots$, we shall mean a bijection $\N \to \N$, which is the identity off a finite set. 
The group of permutations will be denoted by $\Sigma_\infty$. The symmetric group $\Sigma_n$ is identified with the subgroup of $\Sigma_\infty$
consisting of permutations $w$ such that $w(i)=i$ for $i>n$.
We set
$$
A := \mathbb{Z}[x_1,x_2,\ldots]\,.
$$
We define a linear operator $\partial_i: A \to A$ as follows:
$$
\partial_i(f):=
\frac{f(x_1,\ldots ,x_i ,x_{i+1},\ldots ) - f(x_1,\ldots ,x_{i+1},x_i,\ldots )}
{x_i - x_{i+1}}\,.
$$
These are classical Newton's divided differences. For more on divided differences, see \cite{LS1}, \cite{LS2}, \cite{Mcd}, \cite{Lbook}.

For a simple reflection 
$
s_i =1,\ldots ,i-1,i+1,i,i+2,\ldots
$,
we put $\partial_{s_i}:= \partial_i$.

\begin{lemma} Suppose that  $w = s_1  \cdots s_k =
t_1  \cdots  t_k$  are two reduced words for a permutation $w$.
Then we have
$
\partial_{s_1} \circ \cdots  \circ \partial_{s_k}
= \partial_{t_1}\circ \cdots  \circ \partial_{t_k}
$.
\end{lemma}
(See \cite{BGG} and \cite{D}.)

Thus for any $w\in \Sigma_\infty$, we can define  $\partial_w$
as $\partial _{s_1} \circ \cdots  \circ \partial_{s_k}$
independently of a reduced word of $w$.

Let $w\in S_\infty$ and let $n$ be a natural number such that $w(k) = k$ for $k>n$.
We define the {\it Schubert polynomial} of Lascoux and Sch\"utzenberger \cite{LS1}
(1982) associated with a permutation $w$, by setting
$$
{\mathfrak S}_w := \partial_{w^{-1} w_0}(x^{n-1}_1 x^{n-2}_2 \cdots
x^1_{n-1} x^0_{n})\,,
$$
where $w_0$ is the permutation $(n,n-1,\ldots ,2,1),n+1,n+2,\ldots $. Observe that this definition
does not depend on the choice of $n$, because
$$
\partial_n\circ \cdots \circ \partial_2\circ \partial_1(x_1^nx_2^{n-1}\cdots x_n)= x_1^{n-1}x_2^{n-2}\cdots x_{n-1}\,.
$$

\begin{figure}
\includegraphics[width=\textwidth]{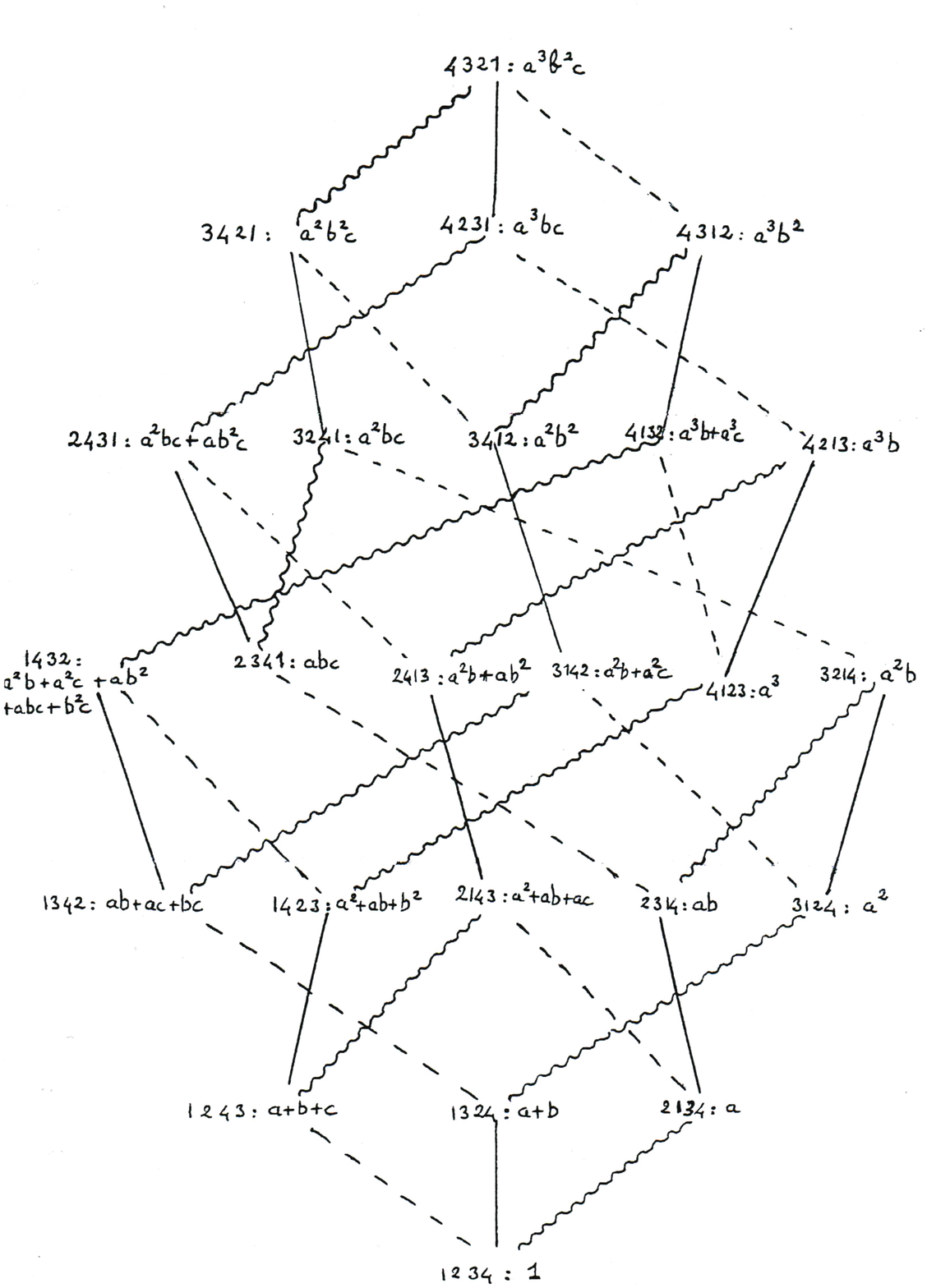}
\end{figure}

\noindent
In the above picture of Schubert polynomials for the symmetric group $\Sigma_4$, we have: $a=x_1, b=x_2, c=x_3$; 
to a permutation there is attached its Schubert polynomial; 
the line $\sim\sim\sim$ means $\partial_1$, the line $\textendash\textendash\textendash$
means $\partial_2$, and the line $- - -$ means $\partial_3$.
(The author got this picture from Lascoux in February 1982 together with a preliminary version of \cite{LS1}.)

\smallskip

Also, we define the $k$th {\it inversion set} of a permutation $w$ as follows:
$$
I_k(w) := \left\{l : l >k,\, w(k) > w(l) \right\} \ \ \ k=1,2,\ldots\,.
$$

By the {\it code} of $w$ (notation: $c(w)$), we shall mean the sequence $(i_1,i_2,\ldots)$, where
$
i_k = |I_k (w)|$, $k=1,2,\ldots
$.
A code determines the corresponding permutation in a unique way.

\begin{example} \rm The code $c(4,2,1,6,3,5,7,8,\ldots)$ is equal
to $(3,1,0,2,0,\ldots)$.
\end{example}

Let us mention the following properties of Schubert polynomials:
\begin{itemize}
\item A Schubert polynomial ${\mathfrak S}_w$ is symmetric in $x_k$ i
$x_{k+1}$ if and only if  $w(k) < w(k+1)$
(or equivalently if $i_k \le  i_{k+1})$;
\item If $w(1) < w(2) < \cdots  < w(k) > w(k+1)
< w(k+2) < \cdots$
(or equivalently $i_1 \le  i_2 \le  \cdots  \le  i_k$,  $0 = i_{k+1} =
i_{k+2}= \cdots)$,
then  $\mathfrak S_w$  is equal to the Schur function $s_{i_k,\ldots ,i_2, i_1}(x_1 ,\ldots ,x_k)$.
\item If  $i_1 \ge i_2 \ge  \cdots$ \,, then $\mathfrak S_w = x_1^{i_1} x_2^{i_2}\cdots$ \
is a monomial.
\end{itemize}

\begin{example} \rm The Schubert polynomials of degree 1 are  ${\frak S}_{s_i}=x_1+\cdots +x_i$, $i=1,2.\ldots$.
\end{example}

For these and more properties of Schubert polynomials, the reader may consult \cite{Mcd} and \cite{Lbook}.

Note that the Schubert polynomials indexed by $\Sigma_n$ do not generate additively $\Z[x_1,\ldots,x_n]$. We have (see \cite{Mcd})
$$
\sum_{w\in \Sigma_n} \Z {\frak S}_w = \bigoplus_{0\le \alpha_i\le n-i} \Z x_1^{\alpha_1} x_2^{\alpha_2} \cdots x_n^{\alpha_n}\,.
$$

\smallskip

If the sets $I_k (w)$ form a chain (with respect to the inclusion),
then the permutation $w$ is called {\it vexillary}. Equivalently, there are no $i < j < k < l$ 
with $w(j)<w(i)<w(l)<w(k)$ (see e.g. \cite{Mcd}).

\smallskip

The following result stems from \cite{LS1} and \cite{Wachs}: 

\begin{theorem}[\rm Lascoux-Sch\"utzenberger, Wachs]
If $w$ is a vexillary permutation with code $(i_1,i_2,\dots,i_n>0,0\ldots)$,
then
$$
\mathfrak S_w=s_{(i_1 ,\ldots ,i_n )^\ge}\left( \min I_1(w)-1 , \ldots  ,
\min I_n(w)-1\right)^\le\,.
$$
\end{theorem}
(By $(-)^{\ge}$ \ and \ $(-)^{\le}$, we mean the respectively ordered sets.)

Here, for two sequences of natural numbers
$$
i_1\ge \cdots \ge i_k \ \ \  \hbox{and} \ \ \  0 < b_1 \le \cdots \le  b_k\,,
$$
$$
s_{i_1,\ldots ,i_k}(b_1,\ldots ,b_k):=
\det \left(
s_{i_p-p+q}(x_1,\ldots ,x_{b_p})\right)_{1\le p,q\le k}
$$
is a {\it flagged Schur function} (see \cite{Mcd}).

\section{Functors asked by Lascoux}\label{las}

Let $R$ be a commutative $\mathbb{Q}$-algebra, and  $E_{\bullet} : E_1 \subset E_2 \subset \cdots$
a flag of $R$-modules. Suppose that $\mathcal I= [i_{k,l}]$,
$k,l =1,2,\ldots$, is a matrix of $0$'s and $1$'s such that
\begin{itemize}

\item $i_{k,l}= 0$ for $k \ge l$;

\item $\sum_l i_{k,l}$ is finite for any $k$;

\item $\mathcal I$ has a finite number of nonzero rows.
\end{itemize}

Such a matrix $\cal I$ is called a {\it shape}:
$$
\begin{matrix}
\scriptstyle0& \scriptstyle0& \scriptstyle0& \scriptstyle0& \scriptstyle0& \scriptstyle0& \scriptstyle1& \scriptstyle0& \ldots\\[-5pt]
\scriptstyle0& \scriptstyle0& \scriptstyle1& \scriptstyle0& \scriptstyle0& \scriptstyle0& \scriptstyle0& \scriptstyle0& \ldots\\[-5pt]
\scriptstyle0& \scriptstyle0& \scriptstyle0& \scriptstyle1& \scriptstyle0& \scriptstyle1& \scriptstyle0& \scriptstyle0& \ldots\\[-5pt]
\scriptstyle0& \scriptstyle0& \scriptstyle0& \scriptstyle0& \scriptstyle1& \scriptstyle0& \scriptstyle1& \scriptstyle0& \ldots\\[-5pt]
\scriptstyle0& \scriptstyle0& \scriptstyle0& \scriptstyle0& \scriptstyle0& \scriptstyle1& \scriptstyle0& \scriptstyle0& \ldots\\[-5pt]
\scriptstyle0& \scriptstyle0& \scriptstyle0& \scriptstyle0& \scriptstyle0& \scriptstyle0& \scriptstyle0& \scriptstyle0& \ldots\\[-5pt]
 & & & \hbox to0pt{$\scriptstyle\ldots\ldots$}
\end{matrix} \ \
= \ \ \begin{matrix}
\scriptstyle 0&\scriptstyle 0&\scriptstyle 0&\scriptstyle 0&\scriptstyle \times&\scriptstyle 0\\[-5pt]
 &\scriptstyle \times&\scriptstyle 0&\scriptstyle 0&\scriptstyle 0&\scriptstyle 0\\[-5pt]
 & & \scriptstyle\times& \scriptstyle0&\scriptstyle \times&\scriptstyle 0\\[-5pt]
 & & &\scriptstyle \times&\scriptstyle 0&\scriptstyle \times\\[-5pt]
 & & & &\scriptstyle \times&\scriptstyle 0\\[-5pt]
 & & & & &\scriptstyle 0
\end{matrix}
$$

Let  $i_k := \sum_{l=1}^\infty i_{k,l}$\,, \ \ \
$\widetilde i_l := \sum_{k=1}^\infty i_{k,l}$\,.
We define $S_{\cal I} (E_\bullet)$  as the image of the following composition:
\begin{equation}\label{Phi}
\Phi_{\cal I}(E_\bullet) : \bigotimes_k S_{i_k}(E_k)
\stackrel{\Delta_S}{\longrightarrow} \bigotimes_k \bigotimes_l
S_{i_{k,l}} (E_k)  \stackrel{m_\wedge}{\longrightarrow }
\bigotimes_l \bigwedge\nolimits^{\widetilde i_l} E_l\,,
\end{equation}
where $\Delta_S$ is the diagonalization in the symmetric algebra and $m_\wedge$ is the multiplication
in the exterior algebra.

For example, if all the 1's occur in $\cal I$ only in the $k$th row, then we get $S_p(E_k)$, where $p$
is the number of 1's. Also, if all the 1's occur in $\cal I$ in one column in consecutive $p$ rows,
then we get $\wedge^p(E_k)$, where $k$ is the number of the lowest row with 1.

Note that $S_{\cal I}(-)$ is a {\it functor}: if $E_\bullet$ and $F_\bullet$ are flags of
$R$-modules, 
$$
f: \bigcup\limits_n E_n \to \bigcup\limits_n F_n
$$ 
is an $R$-homomorphism such that  $f(E_n) \subset F_n$, then $f$ induces
an $R$-homomorphism $S_{\cal I}(E_\bullet) \to  S_{\cal I}(F_\bullet)$.

\begin{remark}\label{vb} \rm \ Note that for a flag $E_{\bullet} : E_1 \subset E_2 \subset \cdots$
of vector bundles on a variety, the same construction leads to a vector bundle $S_{\cal I} (E_\bullet)$ on this variety.
Accordingly, we shall call it a {\it KP bundle}.
\end{remark}

Let $w\in \Sigma_\infty$. By the {\it shape} of $w$ we mean the matrix:
$$
{\cal I}_{w} = \left[i_{k,l}\right]:= \left[\chi_k(l )\right],\
k,l = 1,2,\ldots
$$
where $\chi_k$ is the characteristic function of $I_k(w)$.

\begin{example} \rm
For $w  = 4,2,1,6,3,5,7,8,\ldots $\,, the shape  ${\cal I}_w$
is equal to
$$\scriptstyle{
\begin{matrix}
 \times& \times& 0& \times& 0\cr
 & \times& 0& 0& 0\cr
 & & 0& 0& 0\cr
 & & & \times& \times\cr
 & & & & 0\cr
\end{matrix}}
$$
\end{example}

We define a module $S_w(E_\bullet)$, associated with a
permutation $w$ and a flag  $E_\bullet$ as $S_{{\cal I}_w}(E_\bullet)$; this leads
to a functor $S_w(-)$. These modules were defined by Kra\'skiewicz
and the author in \cite{KP0} and \cite{KP}. Accordingly, we have KP bundles $S_w(E_\bullet)$ associated
with a flag $E_\bullet$ of vector bundles. 

\begin{example} \rm Invoking the previous example, we see that $S_{42163578\ldots}(E_\bullet)$
is the image of the map
$$
\Phi_{{\cal I}_{42163578...}}: S_3(E_1)\otimes S_1(E_2)\otimes S_2(E_4) \to E_1\otimes \wedge^2(E_2)\otimes \wedge^2(E_4)\otimes E_5\,.
$$
\end{example}

\smallskip

From now on, let $E_\bullet$ be a  flag of $K$-vector spaces with
$
\dim E_i=i\,.
$
Let $B$ be the Borel group of linear endomorphisms of
$E := \bigcup  E_i$, which preserve $E_\bullet$.
The modules used in the definition of $S_{\cal I}(E_\bullet)$ are $\mathbb Z[B]$-modules,
and the maps are homomorphisms of $\mathbb Z[B]$-modules.
Let  $\{u_i : i=1,2,\ldots  \}$ be a basis of $E$ such that
$u_1,u_2,\ldots ,u_k$ span $E_k$. Then the module $S_{\cal I}(E_\bullet)$ as a cyclic $\mathbb Z[B]$-submodule
in $\bigotimes_l \bigwedge^{\widetilde i_l}E_l$,
is generated by the element
$$
u_{\cal I}:=  \otimes_l  u_{k_{1,l}} \wedge u_{k_{2,l}} \wedge \cdots
\wedge u_{k_{i_l,l}}\,,
$$
where $k_{1,l} < k_{2,l} < \cdots < k_{i_l,l}$  are precisely those indices
for which $i_{k_{r,l},l}= 1$. In particular, $S_w(E_\bullet)$ is a cyclic $\Z[B]$-module generated by $u_w=u_{{\cal I}_w}$.
The modules $S_w(E_\bullet)$ were called by Watanabe {\it Kra\'skiewicz-Pragacz modules}, in short {\it KP modules} (see \cite[v3]{W1},
\cite{W2}, \cite{W3}).

\begin{example} \rm The KP module $S_{52163478\ldots}(E_\bullet)$
is generated over $\Z[B]$ by the element
$$
u_{42163578\ldots}=u_1\otimes u_1\wedge u_2 \otimes u_1\wedge  u_4 \otimes u_4\,.
$$
\end{example}

KP modules give a substantial generalization of Schur modules discussed in Sect. \ref{schur}.
Note that any Schur module $V_\lambda(E_m)$, where $\lambda=(\lambda_1,\ldots,\lambda_k>0)$ and $k\le m$, 
can be realized as $S_w(E_\bullet)$ for some $w$.
We claim that the permutation $w$ with the code 
$(0^{m-k},\lambda_k,\ldots,\lambda_1,0,\ldots)$ determines the desired KP module. Take for example 
$m=4$, $\lambda=(4,3,1)$. The permutation with the code $(0,1,3,4,0,\ldots)$ is  $1,3,6,8,2,4,5,7,9,\ldots$,
and has shape
$$\scriptstyle{
\begin{matrix}
 0 & 0 & 0 & 0 & 0 & 0 & 0\cr
 & 0 & 0 & \times & 0 & 0 & 0 \cr
 & & 0 & \times & \times & \times & 0 \cr
 & & & \times & \times & \times & \times \cr
 & & & & 0 & 0 & 0 \cr
 & & & & & 0 & 0 \cr
 & & & & & & 0\cr
\end{matrix}}
$$
Its KP module is given by the image of the map
$$
\Phi_{{\cal I}_{136824579,\ldots}}: E_2\otimes S_3(E_3)\otimes S_r(E_4) \to \wedge^3(E_4)\otimes \wedge^2(E_4)\otimes \wedge^2(E_4)\otimes E_4\,.
$$
We invoke at this point a standard basis theorem for Schur modules (see \cite{ABW}, \cite{T}),
which allows us to replace the flag $E_1\subset \cdots \subset E_m$ by $E_m=\cdots = E_m$ without
change of the image of $\Phi_{{\cal I}_w}$. We obtain that the image of $\Phi_{{\cal I}_{136824579,\ldots}}$ is
equal to the image of the map (\ref{sw})
$$
E_4\otimes S_3(E_4)\otimes S_r(E_4) \to \wedge^3(E_4)\otimes \wedge^2(E_4)\otimes \wedge^2(E_4)\otimes E_4\\,
$$
defining $V_{(4,3,1)}(E_4)$.

\smallskip

We shall now study the character of $S_w(E_\bullet)$. In Section \ref{ample}, we shall show an application of the following theorem
to algebraic geometry. In fact, this theorem will tell us what are the Chern roots of a KP bundle as functions
of the Chern roots of the original bundle.

Consider the maximal torus (\ref{T})  $T\subset B$ consisting of diagonal matrices with $x_1,x_2,\ldots$ on the diagonal,
with respect to the basis $\{u_i: i=1,2,\ldots \}$.

\begin{theorem}[\rm Kra\'skiewicz-Pragacz]\label{kp} Let $w\in \Sigma_\infty$. 
The trace of the action of $T$ on
$S_{w}(E_\bullet)$ is equal to the Schubert polynomial $\mathfrak S_w$.
\end{theorem}
About the proof: we study the multiplicative properties of $S_w(E_\bullet)$'s, comparing them with those of the ${\frak S}_w$'s.

Let $t_{p,q}$ be the permutation:
$$
1,\ldots,p-1,q,p+1,\ldots,q-1,p,q+1,... \,.
$$

We now record the following formula for multiplication by $\mathfrak S_{s_k}$:

\begin{theorem}[\rm Monk] Let $w\in \Sigma_\infty$. We have
$$
\mathfrak S_w  \cdot  (x_1 + \cdots  + x_k )  =  \sum \mathfrak S_{w \circ t_{p,q}}\,,
$$
where the sum is over $p,q$ such that $p\le k$, $q>k$  and
$l(w \circ t_{p,q}) =l(w ) + 1$.
\end{theorem}
(See \cite{M}, \cite{LS2}, \cite{Mcd}.)

\begin{example} \rm We have
$$
\mathfrak S_{135246\ldots} \cdot (x_1 +x_2) = \mathfrak S_{235146\ldots}
+ \mathfrak S_{153246\ldots} + \mathfrak S_{145236\ldots}\,.
$$
\end{example}


In fact, we shall use more efficiently the following result of Lascoux and Sch\"utzenberger (see \cite{LS3}).
\begin{theorem}[\rm Transition formula] Let $w\in \Sigma_\infty$.
Suppose that $(j,s)$ is a pair of positive integers such that

\noindent
{\rm 1)} $j<s$  and  $w(j) > w(s)$,

\noindent
{\rm 2)} for any $i\in ]j,s[$,  $w(i) \notin [w(s),w(j)]$

\noindent
{\rm 3)} for any  $r>j$, if $w(s) < w(r)$ then there exists
$i\in ]j,r[$  such that  $w(i) \in  [w(s),w(r)]$.

\smallskip

Then 
$$
\mathfrak S_w = \mathfrak S_v \cdot x_j  + \sum_{p=1}^m \mathfrak S_{v_p}\,,
$$
where $v  = w \circ t_{j,s}$, $v_p = w \circ t_{j,s} \circ t_{k_p,j}$\,,
$p=1,\ldots,m$, say. Here the numbers $k_p$ are such that

\noindent
{\rm 4}) $k_p < j$ and $w(k_p ) < w(s)$,

\noindent
{\rm 5)} if $i\in  ]k_p ,j[$ then $w(i) \notin [w(k_p ),w(s)]$.
\end{theorem}
Note that if $(j,s)$ is the maximal pair (in the lexicographical order) satisfying 1),
then conditions 2)-3) are also fulfilled. A transition corresponding to this pair
will be called {\it maximal}. In particular, for any nonntrivial permutation,
there is at least one transition.

\begin{example} \rm For the permutation $521863479\ldots$,
$$
{\mathfrak S}_{521843679\ldots}
\cdot x_5+{\mathfrak S}_{524813679\ldots}+{\mathfrak S}_{541823679\ldots}
$$
is the maximal transition, and other transitions are
$$
{\mathfrak S}_{521763489\ldots}\cdot x_4+{\mathfrak S}_{527163489\ldots}+{\mathfrak S}_{571263489\ldots}+ {\mathfrak S}_{721563489\ldots}
$$
and 
$${\mathfrak S}_{512864379\ldots} \cdot  x_2\,.
$$
\end{example}

We prove that for the maximal transition for $w$, there exists a filtration
of $\Z[B]$-modules
$$
0  =  \mathcal F_0 \subset \mathcal F_1  \subset \cdots  \subset
\mathcal F_m \subset  \mathcal F  = S_w(E_\bullet)
$$
together with isomorphisms
$$
\mathcal F/\mathcal F_m \simeq S_v (E_\bullet) \otimes E_j /E_{j-1}
 \ \ \ \ \ \ \hbox{and} \ \ \ \ \ \
\mathcal F_p /\mathcal F_{p-1} \simeq
S_{v_p}(E_\bullet)
$$
where $p=1,\ldots ,m$. This implies an isomorphism of $T$-modules
$$
S_w(E_\bullet)\simeq S_v(E_\bullet)\otimes E_j/E_{j-1} \oplus \bigoplus_{p=1}^m S_{v_p}(E_\bullet)\,.
$$
By comparing this with the transition formula for Schubert polynomials, the assertion of Theorem \ref{kp}
follows by a suitable induction.\qed

\smallskip

\noindent
(For details see \cite[Sect. 4]{KP}.)

\bigskip

There exist {\it flagged Schur modules} $S_\la(-)$, associated with suitable shapes (see \cite[p. 1330]{KP}).
Suppose that $E_\bullet$ is a flag of free $R$-modules with $E_1=R$ and such that the $i$th inclusion in the
flag is given by $E_i\hookrightarrow E_i\oplus R \simeq E_{i+1}$. We record

\begin{theorem}[\rm Kra\'skiewicz-Pragacz]
If $w$ is a vexillary permutation with code
$(i_1,i_2,\dots,i_n>0,0\ldots)$,
then
$$
S_w(E_\bullet)=S_{(i_1 ,\ldots ,i_n )^\ge}\left(E_{\min I_1(w)-1},
\ldots, E_{\min I_n(w)-1}\right)^\le\,.
$$
\end{theorem}

\noindent
(See \cite{KP0}, \cite{KP}.)

\smallskip
There is a natural extension of KP modules to complexes  
with interesting applications (see \cite{Sa}).

\section{KP filtrations of weight modules}\label{wat}

All the results of this section (with just a fe(w exceptions) are due to Watanabe in \cite{W1} and \cite{W2}. In \cite{W1}, the author studied the structure of KP modules using the theory of highest weight categories \cite{CPS}.
From the results in \cite{W1}, in particular, one obtains a certain highest weight category whose standard modules are KP modules.
In \cite{W2}, the author investigated the tensor multiplication properties of KP modules.

In this section, we shall use the language and techniques of enveloping algebras, see, e.g., \cite{Di}.

Let $\frak b$ be the Lie algebra of $n\times n$ upper matrices over $K$, $\frak t$ that of diagonal matrices, and $U(\frak b)$ the enveloping algebra of $\frak b$.
Suppose that $M$ is a $U(\frak b)$-module and $\lambda=(\lambda_1,\ldots,\lambda_n)\in {\mathbb Z}^n$,
Denote by 
$$
M_{\lambda}=\{m\in M: hm=\langle \lambda,h \rangle m\}
$$
the {\it weight space} of $\lambda$,  $\langle \lambda,h \rangle =\sum \lambda_i h_i$.
If $M$ is a direct sum of its weight spaces and each weight space has finite dimension, then $M$
is called a {\it weight module}. For a weight module, we set
$$
ch(M):=\sum_\lambda \dim M_\lambda x^\lambda\,,
$$ 
where $x^\lambda=x_1^{\lambda_1}\cdots x_n^{\lambda_n}$.

Let $e_{ij}$ be the matrix with 1 at the $(i,j)$-position and 0 elsewhere.

Let $K_\lambda$ be a one-dimensional $U(\frak b)$-module, where $h$ acts by $\langle \lambda,h \rangle$ and the matrices
$e_{ij}$, where $i<j$, acts by zero.
Any finite-dimensional weight module admits a filtration by these one-dimensional modules.

Fix $n\in \N$. In this section, we shall mainly work with permutations from 
$$
\Sigma^{(n)}: = \{w: w(n+1)<w(n+2)< ...\}\,.
$$
Observe that the codes of permutations in $\Sigma^{(n)}$ are in $\Z_{\ge 0}^n$; in fact, they exhaust $\Z_{\ge 0}^n$.
Moreover, we have
$$
\sum_{w \in \Sigma^{(n)}} \Z \frak S_w = \Z[x_1,\ldots,x_n]\,.
$$

Write $E=\oplus_{1\le i \le n} Ku_i$.
For each $j\in \mathbb N$, let $l_j=l_j(w)$ be the cardinality of the set
$$
\{ i<j: w(i)>w(j)\} = \{i_1<\ldots<i_{l_j} \}\,,
$$
and write
$$
u_w^{(j)}=u_{i_1}\wedge \cdots \wedge u_{i_{l_j}}\in \Lambda^{l_j}(E)\,.
$$
We have
$u_w=u_w^{1}\otimes u_w^{2}\otimes \cdots$ and
$S_w=U(\frak b) u_w$. The weight of $u_w$ is $c(w)$. Observe that Theorem \ref{kp}
can be restated as

\begin{theorem}\label{kp1} For any $w\in \Sigma^{(n)}$, $S_w$ is a weight module
and $ch(S_w)={\frak S}_w$.
\end{theorem}

A natural question arises: What is the annihilator of $u_w$?

\smallskip

Let $w\in \Sigma^{(n)}$. Consider the following assignment:
$$
(1\le i<j\le n) \ \ \longrightarrow \ \  m_{ij}(w)=\#\{ k>j: w(i)<w(k)<w(j)\}\,.
$$
Then $e_{ij}^{m_{ij}+1}$ annihilates $u_w$.
Let $I_w\subset U(\frak b)$ be the ideal generated by $h-\langle c(w), h\rangle$, $h\in \frak t$, and $e_{ij}^{m_{ij}(w)+1}$, $i<j$. Then there exists a surjection
$$
U(\frak b)/I_w \twoheadrightarrow S_w
$$
such that 1 mod $I_w \mapsto u_w$.

\begin{theorem}[\rm Watanabe] This surjection is an isomorphism.
\end{theorem}
(See \cite{W1}, Sect.~4.)

\smallskip
For $\alpha,\beta\in \Z^n$, we shall write $\alpha \pm \beta$ for $(\alpha_1 \pm \beta_1,\ldots,\alpha_n \pm \beta_n)$.

For $\lambda \in {\mathbb Z}_{\ge 0}^n$, we set $S_\lambda:=S_w$, where $c(w)=\lambda$.
For $\lambda \in {\mathbb Z}^n$ take $k$ such that $\lambda+k{\bf 1}\in {\mathbb Z}_{\ge 0}^n$
(${\bf 1}=(1,\ldots, 1)$ $n$ times), and set
$S_\lambda=K_{-k{\bf 1}}\otimes S_{\lambda+k{\bf 1}}$. We shall use a similar notation for Schubert polynomials.

A {\it KP fitration} of a weight module is a sequence
$$
0=M_0\subset M_1\subset \cdots \subset M_r =M
$$
of weight modules such that each subquotient $M_i/M_{i-1}$ is isomorphic to some KP module.

\medskip

One can ask the following questions:

\noindent
1. When a weight module admits a KP filtration?

\smallskip

\noindent
2. Does $S_\lambda \otimes S_\mu$ have a KP filtration?

\medskip

Write $\rho=(n-1,n-2,\ldots,1,0)$. The module $K_\rho$ will play a role of a ``dualizing module''.

Let $\cC$ denote the category of all weight modules. For $\Lambda\subset {\mathbb Z}^n$,
let $\cC_\Lambda$ be the full subcategory of $\cC$
consisting of all weight modules whose weights are in $\Lambda$. If
$|\Lambda|<\infty$ and $\Lambda'=\{\rho - \lambda: \lambda \in \Lambda \}$, then the map
$M \mapsto M^*\otimes K_{\rho}$ yields an isomorphism
$C_{\Lambda'}\cong \cC_{\Lambda}^{op}$.

\begin{lemma} For any $\Lambda \subset {\mathbb Z}^n$, $\cC_{\Lambda}$ has enough projectives.
\end{lemma}
(See \cite[Sect.~6]{W1}.)

\smallskip

We  now define some useful orders on $\Sigma_{\infty}$.
For $w,v\in \Sigma_{\infty}$, $w\le_{lex} v$ if $w=v$ or there exists $i>0$ such that $w(j)=v(j)$ for $j<i$ and
$w(i)<v(i)$.

For $\lambda\in {\mathbb Z}^n$, define $|\lambda|=\sum \lambda_i$. If $\lambda=c(w), \mu=c(v)$, we write $\lambda\ge \mu$
if $|\lambda|=|\mu|$ and $w^{-1}\le_{lex} v^{-1}$.
For general $\lambda, \mu \in {\mathbb Z}^n$ take $k$ such that
$\lambda+k{\bf 1}, \mu+k{\bf 1}\in {\mathbb Z}_{\ge 0}^n$, and define $\lambda\ge \mu$ iff
$\lambda+k{\bf 1}\ge \mu+k{\bf 1}$.

For $\lambda \in {\mathbb Z}^n$, set $\cC_{\le \lambda}:=\cC_{\{\nu: \nu\le \lambda\}}$.

\begin{proposition} For $\lambda\in {\mathbb Z}^n$, the modules $S_\lambda$ and $S_{\rho - \lambda}^*\otimes K_\rho$ are in $\cC_{\le \lambda}$.
Moreover, $S_\lambda$ is projective and $S_{\rho - \lambda}^*\otimes K_\rho$ is injective.
\end{proposition}
(See \cite[Sect.~6]{W1}.)

\smallskip

For the definitions and properties of $\Ext$'s, we refer to \cite{Mc}.
All $\Ext$'s will be taken over $U(\frak b)$, in $\cC_{\le \lambda}$.

\begin{theorem}[\rm Watanabe] For $\mu, \nu \le \lambda$, $\Ext^i(S_\mu, S_{\rho-\nu}^*\otimes K_\rho)=0$ if $i\ge 1$.
\end{theorem}
(See \cite[Sect.~7]{W1}. This can be regarded a ``Strong form of Polo's theorem'' \cite[Theorem 3.2.2]{vK}.)

\begin{theorem}[\rm Watanabe]\label{Ext} Let $M\in \cC_{\le \lambda}$. If $\Ext^1(M, S_{\rho-\mu}^*\otimes K_\rho)=0$ for all $\mu\le \lambda$, then $M$ has a KP filtration such that each of its subquotients is isomorphic to some
$S_\nu$ with $\nu \le \lambda$.
\end{theorem}
(See \cite[Sect.~8]{W1} and also \cite[Theorem 2.5]{W2}.)

\begin{corollary} \label{Ext1} (1) If $M=M_1\oplus \cdots \oplus M_r$, then $M$ has a KP filtration
iff each $M_i$ does.

(2) If $0\to L\to M\to N\to 0$ is exact and $M, N$ have KP filtrations, then $L$ also does.
\end{corollary}
\proof \rm
Assertion (1) follows from $\Ext^1(M,N)=\oplus \Ext^1(M_i,N)$ for any $N$.
\smallskip

\noindent
Assertion (2) follows from the exactness of the sequence
$$
{\Ext}^1(M,A) \to  {\Ext}^1(L,A) \to {\Ext}^2(N,A)
$$
for any $A$.\qed

\begin{proposition}
Let $w\in \Sigma^{(n)}$, $1\le k\le n-1$. Then $S_w \otimes S_{s_k}$ has a KP filtration.
\end{proposition}
This result was established in \cite[Sect. 5]{KP} for $k=1$ and in (\cite[Sect. 3]{W2}) in general.

\begin{theorem}[\rm Watanabe] $S_w\otimes S_v$ has a KP filtration for any $w,v\in \Sigma^{(n)}$.
\end{theorem}
In order to outline a proof, set $l_i=l_i(w)$ and consider an $U(\frak b)$-module
$$
T_w=\otimes_{2\le i \le n}(\Lambda^{l_i}K^{i-1})\,.
$$
The module $T_w$ is a direct sum component of
$$
\otimes_{2\le i\le n} \bigl(S_{s_{i-1}}\otimes \cdots
\otimes S_{s_{i-1}}\bigr)\,,
$$
where $S_{s_{i-1}}$ appears $l_i$ times.

 \begin{proposition} Suppose that $w\in \Sigma_n$. Then there is an exact sequence
$$
0\to S_w \to T_w \to N\to 0\,,
$$
where $N$ has a filtration whose subquotients are $S_u$ with $u^{-1}>_{lex} w^{-1}$.
\end{proposition}

 Granting this proposition, a proof of the theorem consists of considering the following exact sequence of $U(\frak b)$-modules

  $$0 \ \ \to \ \ S_w\otimes S_v \ \ \ \to \ \  T_w\otimes S_v \ \ \ \to \ \ \ N\otimes S_v \ \ \to \ \ 0\,.$$

  Here, the middle module has a KP filtration by the proposition and the module on the right
	one by induction on lex($w$). Consequently the module on the left has a KP filtration
	by Corollary \ref{Ext1}. \qed
	
	To show the proposition, we need the following 
	\begin{lemma}\label{two} (i) If the coefficient of $x^{c(v)}$ in $\frak S_w$ is nonzero, then we have $v^{-1}\ge_{lex} w^{-1}$.
	
	\noindent
	(ii) If $\Ext^1(S_w, S_u)\ne (0)$, then $u^{-1}<_{lex} w^{-1}$.
	\end{lemma}
	(See \cite[Sect. 6]{W1}.)
	
	We come back to the proposition. Define the integers $m_{wu}$ by 
	$$
	\sum_{u\in \Sigma_n} m_{wu} \frak S_u = \prod_{2\le i \le n} e_{l_i}(x_1,\ldots,x_{i-1})\,,
	$$
	where $e_k$ denotes the elementary symmetric function of degree $k$. Thus $m_{wu}$ is the number
	of times $S_u$ appears as a subquotient of any KP filtration of $T_w$. Let us invoke the following
	Cauchy-type formula:
$$
\prod_{i+j\le n}(x_i+y_j)=\sum_{w\in \Sigma_n} \frak S_w(x) \frak S_{ww_0}(y)
$$
(see \cite{Mcd}, \cite{Lbook}). Consider the bilinear form $\langle, \rangle$ on the free abelian group generated
by Schubert polynomials $\frak S_w$, where $w\in \Sigma_n$, corresponding to this Cauchy identity, i.e., such that 
$\langle \frak S_u, \frak S_{u'w_0}\rangle=\delta_{uu'}$ for any $u, u'\in \Sigma_n$.
We have 
$$
m_{wu}=\langle \frak S_{uw_0}, \prod_{2\le i \le n} e_{l_i}(x_1,\ldots,x_{i-1})\rangle\,.
$$
We use now an additional property of the bilinear form: for any $\alpha, \beta \in \Z^n_{\ge 0}$ with $\alpha_i,\beta_i\le n-i$,
$$
\langle x^{\rho-\alpha}, \prod_{1\le i \le n-1} e_{\beta_i}(x_1,\ldots,x_{n-i})\rangle = \delta_{\alpha, \beta}\,.
$$
Using these two mutually orthogonal bases of monomials and products of elementary symmetric functions, we infer that $m_{wu}$ is the coefficient of
$x_1^{n-1-l_n} x_2^{n-2-l_{n-1}} \cdots$ \ in $\frak S _{uw_0}$.
From this description, is not hard to see that for any $k$, 
$$
n-k-l_{n+1-k}=c(ww_0)_k\,,
$$
and thus the number $m_{wu}$ is equal
to the coefficient of $x^{c(ww_0)}$ in $\frak S_{uw_0}$. By Lemma \ref{two}(i), this coefficient is nonzero only if  
$$
u^{-1}\ge_{lex} w^{-1}\,.
$$
If $u=w$, then $m_{wu}=1$; thus the subquotients of any KP filtration of $T_w$ are the KP modules $S_u$, where $u^{-1}>_{lex} w^{-1}$,
together with $S_w$ which occurs just once. Since, by Lemma \ref{two}(ii),  $\Ext^1(S_w,S_u)=0$ if $u^{-1}>_{lex} w^{-1}$, one can take a filtration
such that $S_w$ occurs as a submodule of $T_w$.\qed

One can say that the proposition gives a materialization of the above mentioned Cauchy-type formula for Schubert polynomials.

\begin{theorem}[\rm Watanabe] Let $\lambda\in {\mathbb Z}^n$ and $M\in \cC_{\le \lambda}$. Then we have
  $$
	ch(M) \le \sum_{\nu \le \lambda}{\dim}_K {\Hom}_{U(\frak b)}(M, S_{\rho-\nu}^*\otimes K_\rho)\frak S_\nu\,.
	$$
	(Here $\sum a_\alpha x^\alpha \le \sum b_\alpha x^\alpha$ means that $a_\alpha \le b_\alpha$ for any $\alpha$.)
  The equality holds if and only if $M$ has a KP filtration with all subquotients isomorphic to $S_\mu$, where $\mu \le \lambda$.
  \end{theorem}
	(See \cite[Sect.~8]{W1}.)
	
	\smallskip

As a byproduct, we get a formula for the coefficient of $\frak S_w$ in $\frak S_u \frak S_v$:
  \begin{corollary} \rm This coefficient is equal to the dimension of
 $$
{\Hom}_{U(\frak b)}(S_u\otimes S_v, S_{w_0w}^*\otimes K_\rho)={\Hom}_{U(\frak b)}
	(S_u\otimes S_v\otimes S_{w_ow},K_\rho)\,.
	$$
	\end{corollary}
   \proof We use $ch$ together with its multiplicativity property, and infer
   $$
	\frak S_u \frak S_v = ch (S_u\otimes S_v) =  \sum_w{\dim}_K {\Hom}_{U(\frak b)}(S_u\otimes S_v, S_{\rho-\lambda}^*\otimes K_\rho)\frak S_w\,.\qed
	$$

   \section{An application of KP filtrations to positivity}\label{plet}

   Let $V_\sigma$ denote the Schur functor associated to a partition $\sigma$ (cf. Sect. \ref{schur}), and let $S_\lambda$ be the KP module
	associated with a sequence $\lambda$.

  \begin{proposition} [\rm Watanabe] The module $V_\sigma (S_\lambda)$ has a KP filtration.
	\end{proposition}
	\proof
  The module $(S_\lambda)^{\otimes k}$ has a KP filtration for any $\lambda$ and any $k$.
  Hence 
	$$
	{\Ext}^1((S_\lambda)^{\otimes k}, S_\nu^*\otimes K_\rho)=0
	$$
	for any $\nu$ by Theorem \ref{Ext}. The module
  $V_\sigma(S_\lambda)$ is a direct sum factor of $(S_\lambda)^{\otimes |\sigma|}$.
  Hence 
	$$
	{\Ext}^1(V_\sigma(S_\lambda), S_\nu^*\otimes K_\rho)=0
	$$
	for any $\nu$, and $V_\sigma(S_\lambda)$ has a KP filtration.\qed

  \begin{corollary}\label{coro} \ If $\frak S_w$ is a sum of monomials
	$x^\alpha + x^\beta + \cdots$, then the following specialization of a Schur function: $s_\sigma(x^\alpha, x^\beta, \ldots)$ is a sum of Schubert 
	polynomials with nonnegative coefficients.
	\end{corollary}
	
	\section{An application of KP bundles to positivity}\label{ample}
	
A good reference for the notions of algebraic geometry needed in this section is \cite{Ha1}. For simplicity, by $X$, we shall denote
a nonsingular projective variety.

Given a vector bundle $E$ on $X$ and a partition $\lambda$, by the {\it Schur polynomial of} $E$, denoted $s_\lambda(E)$, we shall mean 
$s_\lambda(\alpha_1,\ldots,\alpha_n)$, 
where $\alpha_i$ are the Chern roots of $E$. 

Recall that a vector bundle $E$ on $X$ is {\it ample} if for any coherent sheaf $\cal F$ on $X$
and $m>>0$, the sheaf $S_m(E)\otimes {\cal F}$ is generated by its global sections (cf. \cite{Ha}).

Recall that KP bundles were introduced in Remark \ref{vb}.
In \cite[p. 632]{F}, Fulton showed that a Schur polynomial of a KP bundle associated to a filtered bundle $E$ is a nonnegative combination 
of Schubert polynomials of $E$. The method relies on ample vector bundles. In fact, a proper variant of KP bundles should
be used (see below).

We shall say that a weighted homogeneous polynomial $P(c_1,c_2,\ldots)$ of degree $d$, where the variables $c_i$ are of degree $i$,
is {\it numerically positive for ample vector bundles} if for any variety $X$ of dimension $d$, 
$$
\int_X P(c_1(E),c_2(E),\ldots)\cap [X]>0\,.
$$
Note that 
under the identification of $c_i$ with the $i$th elementary symmetric function in $x_1,x_2,\ldots$,
any weighted homogeneous polynomial $P(c_1,c_2,\ldots)$ of degree $d$ is a $\Z$-combination of Schur polynomials $\sum_{|\lambda|=d} a_\lambda s_\lambda$.
In their study of numerically positive polynomials for ample vector bundles, the authors of \cite{FL} showed

\begin{theorem} [\rm Fulton-Lazarsfeld] \ A such nonzero polynomial $\sum_{|\lambda|=d} a_\lambda s_\lambda$ is numerically positive for ample
vector bundles if and only if all the coefficients $a_\lambda$ are nonnegative.
\end{theorem}

Schur functors give rise to Schur bundles $V_\lambda(E)$ associated to vector bundles $E$. In \cite[Cor. 7.2]{P}, the following result was shown

\begin{corollary} [\rm Pragacz] \ In a $\Z$-combination $s_\lambda(V_\mu(E))=\sum_\nu a_\nu s_\nu(E)$, all the coefficients $a_\nu$ are nonnegative.
\end{corollary}
Indeed, assuming that $E$ is ample, combining the theorem and the fact that the Schur bundle of an ample bundle is ample (see \cite{Ha}),
the assertion follows.

\smallskip

To proceed further, we shall need a variant of KP modules associated with sequences of surjections of modules. Suppose that
\begin{equation}\label{ff}
F_0=F \twoheadrightarrow F_1 \twoheadrightarrow F_2 \twoheadrightarrow \cdots 
\end{equation}
is a sequence of surjections of $R$-modules, where $R$ is a commutative ring.
Let ${\cal I}=[i_{k,l}]$ be a shape (see Sect. \ref{las}), $i_k := \sum_{l=1}^\infty i_{k,l}$, and
$\widetilde i_l := \sum_{k=1}^\infty i_{k,l}$\,.
We define $S'_{\cal I} (F_\bullet)$  as the image of the following composition:
\begin{equation}\label{Psi}
\Psi_{\cal I}(F_\bullet) : \bigotimes_k \wedge^{i_k}(F_k)
\stackrel{\Delta_\wedge}{\longrightarrow} \bigotimes_k \bigotimes_l
\wedge^{i_{k,l}} (F_k)  \stackrel{m_S}{\longrightarrow }
\bigotimes_l S_{\widetilde i_l}(F_l)\,,
\end{equation}
where $\Delta_\wedge$ is the diagonalization in the exterior algebra and $m_S$ is the multiplication in the symmetric
algebra. For $w\in \Sigma_\infty$, we define 
$$
S'_w(F_\bullet)=S'_{{\cal I}_w}(F_\bullet)\,.
$$
Suppose that $R$ is $K$, and $F_i$ is spanned by $f_1,f_2,\ldots,f_{n-i}$. 
Consider the maximal torus (\ref{T})  $T\subset B$ consisting of diagonal matrices with $x_1,x_2,\ldots$ on the diagonal,
with respect to the basis $\{f_i: i=1,2,\ldots \}$.
If $w\in \Sigma_n$, then the character of $S'_w$, i.e. the trace of the action of $T$ on $S'_w(F_\bullet)$, is equal to ${\frak S}_{w_0ww_0}$.
 
\smallskip

We now pass to filtered vector bundles.
By a {\it filtered bundle}, we shall mean a vector bundle $E$ of rank $n$, equipped with a flag of subbundles
$$
0=E_n\subset E_{n-1}\subset \cdots \subset E_0=E\,,
$$
where $\rank(E_i)=n-i$ for $i=0,1,\ldots,n$.

For a polynomial $P$ of degree $d$, and a filtered vector bundle $E$ on an algebraic variety $X$, 
by substituting $c_1(E_0/E_1)$ for $x_1$, $c_1(E_1/E_2)$ for $x_2$,..., $c_1(E_{i-1}/E_i)$ for $x_i$,...,
we get a class denoted $P(E_\bullet)$. 
Such a polynomial $P$ will be called {\it numerically positive for filtered ample vector bundles}, if for any filtered ample vector bundle on any 
$d$-dimensional variety $X$, 
$$
\int_X P(E_\bullet)\cap [X]>0\,.
$$
Write $P$ as a $\Z$-combination of Schubert polynomials $P=\sum a_w {\frak S}_w$ with the unique coefficients $a_w$.
We then have (see \cite{F})

\begin{theorem} [\rm Fulton] A nonzero such a polynomial $P=\sum a_w {\frak S}_w$ 
is numerically positive for filtered ample bundles if and only if all coefficients $a_w$ are nonnegative.
\end{theorem}
In \cite{F}, in fact, a more general result for $\bf r$-filtered ample vector bundles is proved.

If (\ref{ff}) is a sequence of vector bundles, then our construction gives a vector bundle $S'_w(F_\bullet)$. We record

\begin{lemma} \ If $F$ is an ample vector bundle, then the KP bundle $S'_w(F_\bullet)$ is ample.
\end{lemma}
\proof
The bundle $S'_w(F_\bullet)$ is the image of the map (\ref{Psi}), thus it is a quotient of a tensor product  
of exterior powers of the bundles $F_k$ which are quotients of $F$. The assertion follows from the facts (see \cite{Ha})
that
a quotient of an ample bundle is ample, an exterior power of an ample vector bundle is ample, and the tensor
product of ample bundles is ample.\qed

\smallskip

Consider now a sequence (\ref{ff}) of vector bundles
of ranks $n, n-1,\ldots,1$. Let $x_i=c_1(\Ker(F_{n-i}\to F_{n-i+1}))$. Write $s_\lambda(S'_w(F_\bullet))$ as a $\Z$-combination
\begin{equation}\label{ful}
s_\lambda(S'_w(F_\bullet))=\sum_v a_v {\frak S}_v(x_1,x_2,\ldots)\,.
\end{equation}
Then assuming that $F$ is ample, and observing that $F$ is filtered by the kernels
of the successive quotients,
the following result holds true by the theorem and the lemma.

\begin{corollary} [\rm Fulton] The coefficients $a_v$ in (\ref{ful}) are nonnegative.
\end{corollary}

Let $E$ be a filtered vector bundle with Chern roots $x_1,x_2,\ldots$ and $w$ be a permutation. Then using the notation from the introduction,
the Chern roots of $S_w(E_\bullet)$ are expressions of the form $l(x^\alpha)$, where $x^\alpha$ are monomials of 
${\frak S}_w$. This follows from the character formula for KP modules (cf. Theorem \ref{kp}). Therefore we can restate the corollary in the following way.

\begin{corollary} 
A Schur function specialized with the expressions $l(x^{\alpha})$ associated to the monomials $x^{\alpha}$ of a Schubert polynomial 
is a nonnegative combination of Schubert polynomials in $x_1, x_2,\ldots$. 
\end{corollary}

\medskip

We conclude with some notes.

\begin{note} \rm In a recent preprint \cite{W3}, Watanabe shows that the highest weight category
from \cite{W1} is self Ringel-dual and that the tensor product operation on $U(\frak b)$-modules is compatible
with Ringel duality functor. Moreover, from the paper it follows an interesting corollary about Ext groups
between KP modules: $\Ext^i(S_w, S_v)\simeq \Ext^i(S_{w_0vw_0}, S_{w_0ww_0})$.
\end{note}

\begin{note} \rm It would be interesting to find further applications of KP modules, and construct analogues of KP modules for other types than (A).
\end{note}

\begin{note} \rm In \cite[Sect. 6]{KP}, the reader can find a discussion of some further developments related to KP modules.
\end{note}
	
\noindent
{\bf Acknowledgments} 
The author is greatly indebted to Witold Kra\'skiewicz for collaboration on \cite{KP0, KP}, and for pointing out a defect
in a former version of the present paper.
His sincere thanks go to Masaki Watanabe for helpful discussions, and for pointing out an error 
of exposition in the previous version of this paper.
Finally, the author thanks the organizers of the conference ``IMPANGA 15'' for their devoted work.
\vskip7pt
\noindent
{\bf Note added in proof} \ In a very recent paper \cite{W4}, Watanabe gave explicit KP filtrations materializing
Pieri-type formulas for Schubert polynomials.

\end{document}